\documentclass[12pt,draft]{article}
\usepackage[english]{babel}
\usepackage{amsfonts,amssymb,amsmath}
\newtheorem{definition}{Definition}[section]
\newtheorem{proposition}{Proposition}[section]
\newtheorem{theorem}{Theorem}[section]
\newtheorem{corollary}{Corollary}[section]
\newtheorem{example}{Example}[section]
\newtheorem{remark}{Remark}[section]
\def\halmos{\hfill\rule{6pt}{6pt}}

\begin{document}
\begin{center}
{\Large \bf A new algorithm for computing Pade
approximants}
\end{center}

\begin{center}
{\large Adukov V.M., Ibryaeva O.L.*}

* oli@6v6power.ru (corresponding author)

{\it South Ural State University, Lenin avenue 76, Chelyabinsk,
Russia}
\end{center}

In paper a new definition of {\it reduced Pade approximant} and
algorithm for its computing is proposed. Our approach is based on
the investigation of the kernel structure of the Toeplitz matrix. It
is shown  that the reduced Pade approximant always has nice
properties which classical Pade approximant possesses only in the
normal case. The new algorithm allows us to avoid Froissart doublets
appearance induced by computer roundoff in the non-normal Pade
table.

{\it Pade approximant, Toeplitz matrix, Pade -- Laplace method,
Froissart doublets.}

\section{Introduction}

Pade approximants are locally the best rational approximations
 that can easily be constructed from the
coefficients of a given power series. They are closely related to
continued fractions, orthogonal polynomials and Gaussian quadrature
methods \cite{Brezinski}. They have been widely used in various
problems of  mathematics, physics, and engineering due to their
property to effectively solve the problem of analytic continuation
of the series beyond its disc of convergence \cite{Suetin}.

There are many methods available to compute Pade approximants
\cite{Baker}. Some of them are implemented in computer algebra
systems such as {\it Maple} and {\it Mathematica} and their built-in
utilities are frequently used in applied problems. For calculation
with floating point numbers in {\it Maple} an algorithm due to
Geddes \cite{Geddes} based on fraction-free symmetric Gaussian
elimination is used. Recursive algorithms for computing  Pade
approximants are also widely distributed. At first, they allowed to
find Pade approximants in the case of a normal Pade table, but later
some of them were generalized to the non-normal case; see, for
example, \cite{Bultheel}.

It is well known that the Pade table of a rational function  $R(z)$
is always non-normal since it contains an infinite singular block
which elements are identical to $R(z)$. Actually the entries inside
that block usually differ from the rational function $R(z)$ through
the appearance of supplementary common roots in the numerators and
the denominators, but after reducing the common factors, we get the
rational function $R(z)$. In any practical calculations (because of
the computer roundoff and the noise presence in the input data from
which Pade approximants are constructed), the paired roots in the
numerators and the denominators will not be rigorously equal.  This
phenomenon of \flqq pairing\frqq \ of such zeros and poles got the
name of {\it Froissart phenomenon} and the pairs are known as {\it
Froissart doublets} \cite{Gilewicz}. For example, their appearance
is inevitable in signal processing by using Pade -- Laplace method
\cite{Y}. In order to identify doublets, several Pade approximants
(besides the desired one) are usually calculated. It requires
additional coefficients of the Taylor series. It is undesirable,
since the coefficients are often computed numerically and Pade
approximants are known to be very sensitive to errors in the
coefficients.

The {\it main purpose} of the present paper is to propose a new
algorithm for computing Pade approximants. This algorithm finds the
denominator with minimal degree among all the denominators of the
Pade approximant and is based on the results concerning the kernel
structure of Toeplitz matrices \cite{Adukov98}. In singular block
case it allows to avoid the appearance of Froissart doublets induced by computer roundoff.

The paper is organized as follows. In preliminary Section \ref{Prel} some definitions
and examples are given. These examples show that in order to avoid Froissart doublets
appearance we should study the kernel of the Toeplitz matrix which gives us denominators of the Pade
approximant. It is done in Section \ref{parametrization} which contains the main result
on the parametrization of the set of all the denominators of the Pade approximant.
This result allows us to give the definition of the modified Pade approximant
in Section \ref{modified} and establish its properties. In Section \ref{our}
we propose our algorithm for Pade approximant computing and
provide some examples in order to show that it allows us to avoid
Froissart doublets induced by computer roundoff.

\section{Preliminaries}
\label{Prel}
This section contains some definitions and numerical experiments obtained
in {\it Maple} and {\it
Mathematica}. We have chosen
these packages because they are the most widespread symbolic computation systems
used in research and applications.

The examples illustrate how Froissart
doublets appear or do not appear when we deal with the singular block of a Pade
table for a rational function. Note that our new algorithm for computing Pade approximants will be
proposed in Section \ref{our}. We will implement it in {\it Maple}
system and solve some of the examples again in order to show that the new method allows us to avoid  Froissart
doublets induced by computer roundoff.

We start with the classical definition of Pade approximants.

\begin{definition}
\label{FroPade} (Pade -- Frobenius)

Let $f(z)$ be a (formal) power series $f(z) =
\sum\limits_{k=0}^{\infty} c_k z^k, \ c_k\in\Bbb C$. The $(m,n)$
{\it Pade approximant} corresponding to $f(z)$ is the rational
function $f_{m,n}(z)\allowbreak =
\displaystyle\frac{P_{m,n}(z)}{Q_{m,n}(z)}$, where $P_{m,n}(z)$ and
$Q_{m,n}(z)$ are polynomials in $z$ such that:
\begin{enumerate}
\item \ $Q_{m,n}(z) \not\equiv 0, \ \ \deg Q_{m,n}(z) \le n$,  \ \ $\deg P_{m,n}(z) \le m$,
\item \ $f(z) Q_{m,n}(z) - P_{m,n}(z) = r_{m+n+1} z^{m+n+1}+r_{m+n+2} z^{m+n+2}+\ldots. $
\end{enumerate}
\end{definition}

In a similar way, we can define the Pade approximant at the point
$z=a$ for the series $\sum\limits_{k=0}^{\infty} c_k (z-a)^k.$
Throughout the paper, we will use Pade approximants at the point $z=0.$

Obviously, the coefficients  $q_0,\ldots, q_n$ of the denominator
$Q_{m,n}(z)$ can be obtained by solving the next homogeneous system
of linear equations with a Toeplitz $n\times (n+1)$ matrix:
\begin{equation} \label{bassystem}
\left(\begin{array}{llll}
c_{m+1}&c_m&\ldots&c_{m-n+1}\\
c_{m+2}&c_{m+1}&\ldots&c_{m-n+2}\\
\ \vdots&\ \vdots&&\ \vdots\\
c_{m+n}&c_{m+n-1}&\ldots&c_{m} \end{array}\right)
\left(\begin{array}{l}
q_0\\
q_1\\
\, \vdots\\
q_n \end{array}\right)= \left(\begin{array}{l}
0\\
0\\
\, \vdots\\
0 \end{array}\right).
\end{equation}
Here and further we assume $c_k=0$ if $k<0$. From the Condition 2 it
follows that, with the $q_0,\ldots, q_n$  available, the
coefficients $p_0,\ldots,p_m$ of the numerator $P_{m,n}(z)$ can be
obtained from the formula:
\begin{equation} \label{systemnom}
 \left(\begin{array}{l}
p_0\\
p_1\\
\, \vdots\\
p_m \end{array}\right)=\left(\begin{array}{llll}
c_{0}& 0 &\ldots& 0\\
c_{1}&c_{0}&\ldots&0\\
\ \vdots&\ \vdots&&\ \vdots\\
c_{m}&c_{m-1}&\ldots&c_{m-n}
\end{array}\right)\left(\begin{array}{l}
q_0\\
q_1\\
\, \vdots\\
q_n \end{array}\right).
\end{equation}

Thus, the $(m,n)$  Pade approximant $f_{m,n}(z)$ always exists.
Generally the polynomial $Q_{m,n}(z)$ (and hence $P_{m,n}(z)$) is
found not uniquely, since the rank of the matrix of the system
(\ref{bassystem}) can be less than $n$. Nevertheless, it is easy to
show that the rational function $f_{m,n}(z)$ is unique. Usually
these rational functions are arranged in a double-entry table known
as the {\it Pade table}
 $\left\{f_{m,n}(z)\right\}_{m,n=0,1,\ldots}$ corresponding to
$f(z)$. The Pade approximant $f_{m,n}(z)$ occupies the $n$th line
and the $m$th column of the table. It is well known that identical
Pade approximants can occur only in square blocks of the table. If
the Pade table does not contain such blocks, it is said to be {\it
normal}; otherwise it is called {\it non-normal}.

Provided that there is the denominator $Q_{m,n}(z)$ such that
$Q_{m,n}(0) \ne 0$, we can rewrite the Condition 2 as follows
\begin{equation}\label{BPade}
f(z)  - \frac{P_{m,n}(z)}{Q_{m,n}(z)} = R_{m+n+1}
z^{m+n+1}+R_{m+n+2} z^{m+n+2}+\ldots .
\end{equation}

In this case, the rational function $\frac{P_{m,n}(z)}{Q_{m,n}(z)}$
is called a {\it Pade -- Baker approximant.} Note that the Pade --
Baker approximant does not always exist.

Let us now come to the second and the main part of the section where
we consider how Pade approximants can be obtained in such popular
mathematical packages as {\it Maple, Mathematica} and {\it Matlab}.

The built-in utility {\it pade} in {\it Matlab} approximates time
delays by rational linear-time invariant models, i.e. actually it
finds only diagonal Pade approximants for an exponential function by
available for this case explicit formulas \cite{Golub2}. {\it Maple}
and {\it Mathematica} compute Pade approximants in the general case
and in the examples below we will use their functions {\it pade($f$,
$z=a$, $[m,n]$)} (in {\it Maple}) and  {\it PadeApproximant\,$[\,
expr,\,{z, \,a, \,{m,n}}\,]$} (in {\it Mathematica}) which compute
the $(m,n)$ Pade approximant of the function $f(z)$ at the point
$z=a$.

To compare the results with the approximated rational functions we will find zeros and poles of approximants.

\begin{example}\label{noFroissart}
Let us find by {\it Maple} and {\it Mathematica} the diagonal $(2,2)$ Pade
approximant $f_{2,2}$ for $f(z)=\frac{(z+1)(z-2)}{(z+2.1)(z-1)}$ at the point
$z=0$ via the following commands.

In {\it Maple} we have:

$f := \frac{(z+1)\cdot(z-2)}{(z+2.1)\cdot(z-1)}:$

with(numapprox): $p := pade(f, z = 0, [2, 2]);$
\begin{center}
$\frac{0.9523809523+0.4761904764 z-0.4761904759z^2}{0.9999999999-0.5238095237z-0.4761904759 z^2}$
\end{center}

fsolve(numer(p), z, complex);
\vspace{-2mm}
\begin{center}
\footnotesize $-1.000000000, 2.000000001$
\end{center}
\vspace{-2mm}

fsolve(denom(p), z, complex);
\vspace{-2mm}
\begin{center}
\footnotesize $-2.100000001, 1.000000000$
\end{center}
\vspace{-2mm}

In  {\it Mathematica} we have:

In[1]:= \ f := (z + 1)*(z - 2)/((z + 2.1)*(z - 1));

In[2]:=\ P := PadeApproximant[f, \{z, 0, \{2, 2\}\}]; P

Out[2]:= \ $\frac{0.952381 + 0.47619 z - 0.47619 z^2}{1.000000000000000 - 0.52381 z -
 0.47619 z^2}$

In[3]:=\ Roots[Numerator[P] == 0, z]

Out[3]:=\ $z == -1. || z == 2.$

In[4]:=\ Roots[Denominator[P] == 0, z]

Out[4]:=\ $z == -2.1 || z == 1.$

\end{example}

As it can be seen,  the obtained approximant $f_{2,2}$ is identical to the function $f(z)$ since they have the same zeros and poles. Note that in this case the kernel dimension of the matrix of system (\ref{bassystem})
is equal to 1, i.e. denominator of the Pade approximant is unique.

It can be verified that $f_{3,3}$ is again identical to the function $f(z)$ (both in {\it Maple} and in {\it Mathematica}), though the kernel is multidimensional.

However, approximant $f_{4,4}$ does not coincide with $f(z)$ because of an appearance of supplementary roots of the denominator and the numerator.

\begin{example} \label{Fr}
Let us find the  $(4,4)$ Pade
approximant for the function $f(z)=\frac{(z+1)(z-2)}{(z+2.1)(z-1)}$ at the point
$z=0$.

In {\it Maple} we have:

p := pade(f, z = 0, [4, 4]);

\begin{center}
$\frac{0.9523809524+0.9750566894 z-0.2267573697 z^2-0.2494331066 z^3}{1.-
0.7505668935 z^2-0.2494331066 z^3}$
\end{center}

fsolve(numer(p), z, complex);
\vspace{-2mm}
\begin{center}
\footnotesize  ${\bf -1.909090909}, -.9999999999, 2.000000000$
\end{center}
\vspace{-2mm}

fsolve(denom(p), z, complex);
\vspace{-2mm}
\begin{center}
\footnotesize $-2.100000001, {\bf -1.909090908}, 1.000000000$
\end{center}

In {\it Mathematica} we have:

In[1]:= \ P := PadeApproximant[f, \{z, 0,\{4, 4\}\}]; P

Out[1]:=\
$\frac{0.952381 + 1.11851 z - 0.155032 z^2 - 0.321159 z^3 -
 3.33067\times 10^{-16} z^4}{1.000000000000000 + 0.150624 z - 0.829465 z^2 -
 0.321159 z^3 + 0. z^4}$

In[2]:=\ Roots[Numerator[P] == 0, z]

Out[2]:=\ z == $-9.64247\times 10^{14} ||  z =={\bf -1.48273} || z == -1. || z == 2.$

In[3]:= \ Roots[Denominator[P] == 0, z]

Out[3]:=\ $z == -2.1 || z == {\bf -1.48273} || z == 1.$
\end{example}

Thus, in Example \ref{Fr} Froissart doublets induced by computer roundoff appear. Note that doublets
 obtained in {\it Maple} and {\it Mathematica} are different (they are in bold).
The root $ -9.64247\times 10^{14}$ of the numerator (obtained in {\it Mathematica})
is caused by its small leading coefficient $-3.33067\times 10^{-16}$.

Note that for $f_{3,3}$ and $f_{4,4}$ the denominator is not unique
and, as we have seen, {\it Maple} and {\it Mathematica} may choose not the best one.

Let us take another rational function and find
its approximation in its infinite singular block.
As it can be seen below,
{\it} Maple does not produce Froissart doublets, but they appear in {\it Mathematica}.

\begin{example}\label{Fr1}
Let us find the $(2,3)$ Pade approximant for the function
$f(z) = \frac{z+1.01}{(z+2)(z-2.01)}$ at the point $z=0.$

In {\it Maple} we have:

f := $\frac{z+1.01}{(z+2)\cdot(z-2.01)}:$

with(numapprox): p := pade(f, z, [2, 3]);

\begin{center}
$\frac{-0.2512437812-0.2487562190 z}{1.000000000+0.002487561992 z-0.2487562188 z^2}$
\end{center}

fsolve(numer(p), z, complex);
\vspace{-2mm}
\begin{center}
\footnotesize $ -1.010000000$
\end{center}
\vspace{-2mm}

fsolve(denom(p), z, complex);
\vspace{-2mm}
\begin{center}
\footnotesize $-2.000000001, 2.010000000$
\end{center}

In {\it Mathematica} we have:

In[1]:=\ f := (z + 1.01)/((z + 2)*(z - 2.01));

In[2]:=\ P := PadeApproximant[f, \{z, 0, \{2, 3\}\}]; P

Out[2]:=\ $\frac{-0.251244 - 0.122069 z + 0.125433 z^2}{1.000000000000000 -
 0.501752 z - 0.250011 z^2 + 0.125433 z^3}$

In[3]:=\ Roots[Numerator[P] == 0, z]

Out[3]:=\ $z == -1.01 || z == {\bf 1.98318}$

In[4]:=\ Roots[Denominator[P] == 0, z]

Out[4]:=\ $z == -2. || z == {\bf 1.98318} || z == 2.01$
\end{example}

The following conclusions may be drawn from the examples.
When the kernel dimension of the matrix of system (\ref{bassystem})
is equal to 1 (as in Example \ref{noFroissart}), we have one denominator of the Pade approximant and
Froissart doublets can not appear. But they can sometimes appear  when the kernel
dimension is greater than 1.
If an improper denominator was chosen,
then supplementary (also called artificial) roots appear.
These roots may have corresponding pairs among the roots of the
numerator and the cancelation of the common factors reduces the Pade approximant
to the reduced form. But in the presence of computer roundoff,
the paired roots in the numerator
and denominator will not be rigorously equal and
this can cause the appearance of Froissart doublets (as in Examples \ref{Fr}, \ref{Fr1}). Moreover, computer roundoff may produce artificial roots with the great absolute value (as in Example \ref{Fr}), since the vanishing
leading coefficients are small but not equal to zero.

To avoid these problems, we should choose the denominator with the minimal degree.
As will be shown further, this denominator always exists. In order
to find it, we have to study the kernel structure of the Toeplitz matrix
of the system (\ref{bassystem}). It will be done in the following section.

\section{Parametrization of the denominator set}
\label{parametrization}

 In this section we study the structure of the set of all denominators for a $(m,n)$ Pade approximant, i.e. the structure of the kernel    ${\rm ker}\,T_{m+1}$
 of the Toeplitz matrix $$T_{m+1}=\|c_{i-j}\|_{\mbox{\scriptsize${\begin{array}{l}
i=m+1,\ldots,m+n\\
j=0,1,\ldots,n\end{array}}$}}$$ from the system (\ref{bassystem}).

First of all, let us prove that the minimal degree denominator exists and establish its properties.

\begin{proposition}\label{existmin}
There exists the denominator $Q_{m,n}^{0}(z)$ with the minimal
degree among all the denominators $Q_{m,n}(z)$. This denominator is
unique up to a constant factor.

Let $P_{m,n}^{0}(z)$ be the numerator of the Pade approximant corresponding to
$Q_{m,n}^{0}(z)$. Then the polynomials $P_{m,n}^{0}(z)$ and $Q_{m,n}^{0}(z)$
do not have common non-zero roots.
\end{proposition}
{\bf Proof.} Since all the polynomials $Q_{m,n}(z)$ satisfy the condition
 $\deg Q_{m,n}(z)\leq n$, there exists the denominator with the minimal degree $d$.

Suppose that there are two denominators with the degree $d$:
$$
Q_{m,n}^{0}(z)=B_d z^{d}+\ldots+B_0,  \ \
\widetilde{Q}_{m,n}^{0}(z)=\widetilde{B}_d z^{d}+\ldots+\widetilde{B}_0, \ \  B_d\neq0, \ \widetilde{B}_d\neq0.
$$

Let us introduce $Q(z)=\widetilde{B}_d Q_{m,n}^{0}(z)-B_d \widetilde{Q}_{m,n}^{0}(z)$.
This is the denominator of the Pade approximant with the degree less than $d$. It is
possible if and only if $Q(z)\equiv0$. Thus, $Q_{m,n}^{0}(z)$ and
$\widetilde{Q}_{m,n}^{0}(z)$ are linearly dependent. The uniqueness is proved.

Let $P_{m,n}^{0}(z)$ be the numerator corresponding to $Q_{m,n}^{0}(z)$.
Let us suppose that $z_{0}$ is the common non-zero root of $P_{m,n}^{0}(z)$ and $Q_{m,n}^{0}(z)$:
$$
P_{m,n}^{0}(z)=(z-z_{0}) \widetilde{P}_{m,n}^{0}(z), \ \
Q_{m,n}^{0}(z)=(z-z_{0}) \widetilde{Q}_{m,n}^{0}(z).
$$
By Definition~\ref{FroPade} we have:
$$
f(z) \widetilde{Q}_{m,n}^{0}(z)-\widetilde{P}_{m,n}^{0}(z)=
\frac{1}{\left(z-z_{0}\right)}\, O \left(z^{m+n+1}\right)={
O}\left(z^{m+n+1}\right).
$$
\\
Then $\widetilde{Q}_{m,n}^{0}\left(z\right)$ is the denominator with
the degree less than $d$. We have a contradiction since $d$ is the minimal degree.

\halmos

Now we describe the denominator set $\{Q_{m,n}(z)\}$, i.e. the kernel structure of $T_{m+1}$. This result will be used to construct the algorithm for obtaining the minimal degree denominator.

To find any denominator of the Pade approximant, we need the
sequence of the Taylor coefficients  $c_{m-n+1},\ldots,c_{m+n}$ of
$f(z)$. The sequence consists of the entries of matrix $T_{m+1}$.
Let us introduce notions  of indices and essential polynomials for
$c_{m-n+1},\ldots,c_{m+n}$. The notions were given in more general
situation in the paper~\cite{Adukov98}. Here we formulate them
specially for our case.

It is natural to include $T_{m+1}$ in the family of Toeplitz matrices
\begin{equation}
T_{k}=
\begin{pmatrix}
    c_{k}   & c_{k-1} &\ldots & c_{m-n+1}  \\
    c_{k+1} & c_{k}   &\ldots & c_{m-n+2}\\
    \vdots  & \vdots  &\ddots &\vdots  \\
    c_{m+n}   & c_{m+n-1} &\ldots & c_{2m-k+1}
\end{pmatrix}, \ \ \ m-n+1 \leq k \leq m+n,
\end{equation}
which are constructed from the elements of $c_{m-n+1},\ldots,c_{m+n}$.

Consider the sequence of the spaces  ${\ker}\,T_{k}$.
It is more convenient to deal not with  vectors $Q=(q_0, q_1, \ldots, q_{k-m+n-1})^t \in {\ker}\,T_{k}$ but with their generating polynomials $Q(z)=q_{0}+q_{1}\,z+\cdots+q_{\,k-m+n-1}\,z^{k-m+n-1}$. Instead of the spaces $\ker T_k$ we will use the isomorphic spaces ${\cal N}_{k}$ consisting of the generating  polynomials.

To do this, we introduce a linear functional $\sigma$  by the formula:
$$
\sigma\{z^j\}=c_{-j}, \ \ -m-n\leq j\leq -m+n-1.
$$
In the theory of orthogonal polynomials this functional is called the Stieltjes functional.

Denote by ${\cal N}_{k} \ (m-n+1 \leq k \leq m+n)$
the space of polynomials $Q(z)$ with the formal degree
$k-m+n-1$ satisfying the orthogonality conditions:
\begin{equation}  \label{e-kerR}
\sigma \bigl\{z^{-i}Q(z)\bigr\}=0,\ i=k,k+1,\ldots,m+n.
\end{equation}

It is easily seen that  ${\cal N}_k$ is the space of generating polynomials  of vectors in $\ker T_k$. For convenience, we put
${\cal N}_{m-n}=0$ and denote by ${\cal N}_{m+n+1}$ the $(2n+1)$-dimensional space of all polynomials with the formal degree $2n$.

Let $d_{k}$ be the dimension of the space ${\cal N}_{k}$ and
$\Delta_{k}=d_{k}-d_{k-1}\ (m-n+1\leq k \leq m+n+1)$. The following fact is crucial for the further considerations.

\begin{proposition}  \label{t-sker}
For any non-zero sequence $c_{m-n+1},\ldots,c_{m+n}$  the following inequalities
\begin{equation}
0 = \Delta_{m-n+1}\leq \Delta_{m-n+2} \leq
\ldots \leq \Delta_{m+n} \leq \Delta_{m+n+1} = 2 \label{e-rchain}
\end{equation}
are fulfilled.
\end{proposition}
{\bf Proof.} It follows from Definition~\ref{FroPade} that ${\cal N}_k$ è $z{\cal N}_k$ are subspaces of ${\cal
N}_{k+1}\ (m-n\leq k\leq m+n)$ and
\[ {\cal N}_{k} \cap z{\cal N}_{k} = z{\cal N}
_{k-1}. \]
Hence, by the Grassman formula, we have:
\begin{equation}   \label{e-dimsum}
\dim \bigl({\cal N}_{k} + z{\cal N}_{k}\bigr) =
2d_{k} - d_{k-1}.
\end{equation}

Let $h_{k+1}$ be the dimension of any complement  ${\cal
H}_{k+1}$ of the subspace ${\cal N}_{k} + z{\cal
N}_{k}$ in the whole space ${\cal N}_{k+1}$.
From~(\ref{e-dimsum}) we have
\begin{equation} \label{hk}
h_{k+1} = \Delta_{k+1} -
\Delta_{k},
\end{equation}
 i.e. $\Delta_{k+1} \geq
\Delta_{k}$. By definition $d_{m-n}=0$ and $d_{m-n+1}$ is also equal to zero.
Hence, $\Delta_{m-n+1}= 0$. In a similar manner we can prove that $\Delta_{m+n+1} = 2$.
\halmos

It follows from the inequalities (\ref{e-rchain}) that there exist
integers $\mu_{1}\leq \mu_{2}$ such that
\begin{equation}   \label{e-DeltaR}
\begin{array}{ccccccl}
\Delta_{m-n+1}&=&\ldots&=&\Delta_{\mu_{1}}&=&0, \\
\Delta_{\mu_{1}+1}&=&\ldots&=&\Delta_{\mu_{2}}&=&1, \\
\Delta_{\mu_{2}+1}&=&\ldots&=&\Delta_{m+n+1}&=&2.
\end{array}
\end{equation}
If the second row in these relations is absent, we assume $\mu_1=\mu_{2}$. Really in our case $\mu_1<\mu_2$ as will be shown later.

\begin{definition} \label{indices}
The integers $\mu_1,\mu_{2}$ defined in~(\ref{e-DeltaR}) will be called
 the essential indices (briefly, indices) of the sequence $c_{m-n+1},\ldots,c_{m+n}$.
\end{definition}

\begin{proposition} \label{indformula}
Let $\varkappa={\rm rank}\,T_{m}$. Then the indices $\mu_1,\mu_2$
are found by the formulas:
\[ \mu_1=m-n+\varkappa, \ \ \ \mu_2= m+n -\varkappa+1. \]
\end{proposition}
{\bf Proof.}
It follows from the definition of  $\Delta_k$ that $\sum\limits_{k=m-n+1}^{m+n+1}\Delta_k=d_{m+n+1}-d_{m-n}=2n+1. $

On the other hand, from the relations (\ref{e-DeltaR}) we have
\[\sum_{k=m-n+1}^{m+n+1}\Delta_k=1\cdot(\mu_2-\mu_1)+2\cdot(m+n+1-\mu_2).
\]
Hence, $\mu_1+\mu_2=2m+1.$ It means  that $\mu_1\leq m < \mu_2.$

In a similar manner we obtain
$\sum\limits_{k=m-n+1}^{m}\Delta_k=d_{m}=n-\varkappa$  and
$\sum\limits_{k=m-n+1}^{m}\Delta_k=m-\mu_1.$ Thus,
$\mu_1=m-n+\varkappa.$

Since $\mu_1+\mu_2=2m+1$, we get $\mu_2= m+n
-\varkappa+1.$

 \halmos

Now we can describe the structure of the kernels of the matrices $T_k$.
It follows from~(\ref{hk}) and~(\ref{e-DeltaR}) that
$h_{k+1}:=\dim {\cal H}_{k+1} \neq 0$ iff $k=\mu_j\ (j=1,2)$. In that case
$h_{k+1}=1.$  Therefore, for $k\neq \mu_j$
\begin{equation}
{\cal N}_{k+1}={\cal N}_k+z{\cal N}_k,
\label{e-rker1}
\end{equation}
and for $k = \mu_{j}$
\begin{equation}
{\cal N}_{k+1} = \bigl({\cal N}_k+z{\cal
N}_k\bigr) \dot{+} {\cal H}_{k+1}.  \label{e-rker2}
\end{equation}

\begin{definition}
Any polynomial $Q_j(z)$ that forms a basis for one-dimensional
complement ${\cal H}_{\mu_{j}+1}$ will be called the essential
polynomial of the sequence $c_{m-n+1},\ldots,c_{m+n}$ corresponding
to the index $\mu_j,$ $j=1,2.$
\end{definition}

It can be shown (see \cite{Adukov98}, theorem 4.1) that integers
$\mu_1,\mu_2$ are the indices and polynomials  $Q_1(z)\in {\cal
N}_{\mu_1+1}$, $Q_2(z)\in {\cal N}_{\mu_2+1}$ are the essential
polynomials iff
\begin{equation} \label{sigma0}
\sigma_0:=\sigma\{z^{-m-n-1}[Q_2(0)Q_1(z)-Q_1(0)Q_2(z)]\}\ne 0.
\end{equation}

In the following theorem the structure of the kernels of matrices
$T_k\ (m-n+1\leq k\leq m+n)$ is described in terms of the indices
and the essential polynomials.

\begin{theorem}   \label{t-basis}

Let $\mu_1,\mu_{2}$ be the indices and let
$Q_{1}(z),Q_{2}(z)$  be the essential polynomials of the sequence
$c_{m-n+1},\ldots,c_{m+n}$.

Then
\begin{center}
$ {\cal N}_{k}=
\left\{
\begin{array}{c}
0, \\
\{q_{1}(z) Q_{1}(z)\}, \\
\{q_{1}(z) Q_{1}(z)+q_{2}(z) Q_{2}(z)\},
\end{array} \right.
$ $
\begin{array}{c}
m-n+1 \le k \le \mu_{1},  \\
\mu_{1}+1 \le k \le \mu_{2}, \\
\mu_{2}+1 \le k \le m+n,
\end{array}
$
 \end{center}
where $q_{1}(z)$, $q_{2}(z)$ are arbitrary polynomials of the formal degree
 $k-\mu_{j}-1$.
\end{theorem}
{\bf Proof.} From~(\ref{e-DeltaR}) we have $d_k=0$ for
$k\in[m-n+1,\mu_1)$, i.e. ${\cal N}_k=0$.

Let $\mu_{1}+1 \le k \le \mu_{2}$.
It follows from~(\ref{e-rker1}) and (\ref{e-rker2}) that the polynomials
\begin{equation}
\left\{ R_{1}(z),\ zR_{1}(z),\ \ldots,\ z^{k-\mu_{1}-1}R_{1}(z)
\right\}       \label{e-rbas}
\end{equation}
generate the space ${\cal N}_k$. The number of these polynomials is equal to
$k-\mu_1$.

From the definition of $\Delta_j$ and the relations (\ref{e-DeltaR}) we have
$d_k=\sum\limits_{j=m-n+1}^k\Delta_j=k-\mu_1$.

Thus, the number of the polynomials in~(\ref{e-rbas}) is equal to
the dimension of the space ${\cal N}_k$. Therefore, these polynomials form a basis of ${\cal N}_k$.

The case $\mu_{2}+1 \le k \le m+n$ can be considered in a similar manner.

  \halmos

Now we apply the previous theorem to the case  $k=m+1$ and
obtain the main result of this section on the parametrization of the denominator set.

\begin{theorem}
\label{param}
Let $Q_{m,n}(z)$ be an arbitrary denominator of the $(m,n)$ Pade approximant
for the series $f(z)=\sum\limits^{\infty}_{k=0}c_{k}z^{k}$. Form the sequence
$\left\{c_{m-n+1},\ldots, c_{m+n}\right\}$ which is necessary to find
$Q_{m,n}(z)$.

Let $\mu_{1}$ be the first index, let $Q_{1}(z)$ be the first essential polynomial of this sequence.

Then the denominator set is
$$\left\{Q_{m,n}(z)\right\}=\left\{q_{1}(z) {Q_{1}(z)}\right\},$$ where
$q_{1}(z)$ is an arbitrary polynomial with the formal degree $m-\mu_{1}$.
Thus, $Q_{1}(z)$ is the denominator $Q_{m,n}^{0}(z)$ with the minimal degree.
\end{theorem}

\begin{remark} \label{PB}

It follows from $Q_{m,n}(z)=q(z)Q_1(z)$ that
the denominator $Q_{m,n}(z)$ such that $Q_{m,n}(0)\ne
0$ (Baker's condition) exists iff $Q_1(0)\ne 0$.

Hence, the Pade approximant with the minimal degree denominator
is the Pade -- Baker approximant if the latter exists.

\end{remark}

In the following section we use this parametrization in order to modify the definition of the Pade approximants.

\section{Reduced Pade approximant and its properties}
\label{modified}
Definition \ref{FroPade}
ignores the non-uniqueness  (in general case) of a denominator of a Pade approximant.
It does not matter much for the normal Pade table but in the singular case
common factors of the denominator and the nominator can not be cancelled because of roundoff errors. This is one of reasons why Froissart doublets can appear.

The main aim of this section is to modify the classical Definition~\ref{FroPade}.
We would like a Pade approximant to exist always, to have a unique denominator (up to a constant factor) and to coincide with the Pade -- Baker approximant if the latter exists.

Let us add the minimality requirement to Definition \ref{FroPade}.

\begin{definition} \label{modpade}
Let $f(z)$ be a (formal) power series $f(z) =
\sum\limits_{k=0}^{\infty} c_k z^k, \ c_k\in\Bbb C$. The $(m,n)$
{\it Pade approximant} corresponding to $f(z)$ is the rational
function $f_{m,n}(z)\allowbreak =
\displaystyle\frac{P_{m,n}(z)}{Q_{m,n}(z)}$, where $P_{m,n}(z)$ and
$Q_{m,n}(z)$ are polynomials in $z$ such that:
\begin{enumerate}
\item \ $Q_{m,n}(z) \not\equiv 0, \ \ \deg Q_{m,n}(z) \le n$,  \ \ $\deg P_{m,n}(z) \le m$,
\item \ $f(z) Q_{m,n}(z) - P_{m,n}(z) = r_{m+n+1} z^{m+n+1}+r_{m+n+2} z^{m+n+2}+\ldots. $
\item The polynomial $Q_{m,n}(z)$ has the minimal degree among
all polynomials satisfying 1, 2.
\end{enumerate}

\end{definition}

Theorem \ref{param} on parametrization of the denominator set gives the constructive method for finding the denominator with the minimal degree.

The next theorem shows that the Pade approximant has the desired properties.

\begin{theorem} \label{proper}

For any power series $f(z)$ the following statements are fulfilled.

\begin{enumerate}

    \item Any $(m,n)$ Pade approximant exists and is unique.

    \item The denominator $Q_{m,n}(z)$ of the Pade approximant is unique up to a constant factor and $Q_{m,n}(z)$ is the first essential polynomial $Q_{1}(z)$ of the sequence $\left\{c_{m-n+1},\ldots,c_{m+n}\right\}$.

        \item $P_{m,n}(z)$ and $Q_{m,n}(z)$ have not common non-zero roots.

        \item The Pade - Baker approximant exists iff $Q_{1}(0)\neq 0$. The Pade approximant from Definition \ref{modpade} is the Pade -- Baker approximant if the latter exists.

            \item If $Q_{1}(z)$ has the root $z=0$ of order $\delta_{m,n}>0$, then
            \[f(z)-f_{m,n}(z)=A z^{m+n+1-\delta_{m,n}}+B
                z^{m+n+2-\delta_{m,n}}+\ldots, \ \ A \ne 0.\]
                \\ ($\delta_{m,n}$ is called the deficiency index of Pade approximant.)

            \end{enumerate}

\end{theorem}
{\bf Proof.} Statements 1 -- 4 evidently follow from
Proposition~\ref{existmin}, Theorem~\ref{param} and Remark~\ref{PB}.

Prove the last statement of the theorem. Let $z=0$ be the root of order
$\delta_{m,n}$ of $Q_{1}(z)$. Since $Q_1(0)=0$, the number $\sigma_0$  from formula
(\ref{sigma0}) is
\[
\sigma_0=-\sigma\{z^{-m-n-1}Q_1(z)\}Q_2(0) \ne 0.
\]
Here $Q_2(z)$ is the second essential polynomial.

Thus $\sigma\{z^{-m-n-1}Q_1(z)\} \ne 0.$

It is easy to see that the number $\sigma\{z^{-m-n-1}Q_1(z)\}$ is the coefficient at $z^{m+n+1}$ in the power series $f(z)Q_1(z)$. As it is not zero, we have
$$
f(z)Q_1(z) - P_{1}(z) = A_1z^{m+n+1}+\ldots \ , \ \ \ \
A_1=\sigma\{z^{-m-n-1}Q_1(z)\} \ne 0.
$$
Here $P_1(z)$ is the numerator corresponding to the denominator $Q_1(z)$.

After dividing this equation by $Q_1(z)$, we get
$$
f(z) - f_{m,n}(z) = Az^{m+n+1-\delta_{m,n}}+\ldots \ ,
$$
where $A\ne 0$.

It means that the deficiency index of the Pade approximant $f_{m,n}(z)$ coincide with the multiplicity
$\delta_{m,n}$ of the root $z=0$ of the first essential polynomial $Q_1(z)$. \halmos

Since the numerator and the denominator of Pade approximant from Definition \ref{modpade} have not common non-zero roots, we will call it {\it reduced Pade approximant}. Note that the numerator and the denominator of the reduced Pade approximant can not be coprime because of their common zero roots.

\begin{remark} It is impossible to improve Definition \ref{modpade} in such a way that
the denominator and the numerator are always coprime.

It is easy to see that if $Q_{m,n}(z)=z^{\delta}Q_{m,n}^0(z)$, then $P_{m,n}(z)=z^{\delta}P_{m,n}^0(z)$.
Hence  $f_{m,n}(z)=\displaystyle\frac{P_{m,n}^0(z)}{Q_{m,n}^0(z)}$.
However  the polynomials $P_{m,n}^0(z)$, $Q_{m,n}^0(z)$ can not be considered as the numerator and the denominator of $f_{m,n}(z)$, since  the Frobenius Condition 2 does not fulfilled:
$$
f(z)Q_{m,n}^0(z) - P_{m,n}^0(z) = Az^{m+n-\delta+1}+\ldots.
$$
\end{remark}

To determine $\delta$ (the deficiency index of $f_{m,n}(z)$), we should know how many low order coefficients of $Q_1(z)$ are equal to zero. Note that vanishing coefficients of the denominator and the numerator can be non-zero because of roundoff errors.
In particulary, it can cause the appearance of roots of the denominator and/or the numerator with the great absolute value (as in Example \ref{Fr}).

Thus, we face the problem of finding vanishing coefficients of the denominator and the numerator.
The results of Theorem \ref{proper} allow to solve this problem.
This will be done in the following two theorems.

Recall that $\varkappa={\rm rank}\, T_m$, $\mu_1=m-n+\varkappa$, matrix
$T_{\mu_1+1}$ has the size  $(2n-\varkappa) \, \times\, (\varkappa+1)$, and
${\rm rank}\, T_{\mu_1+1}=\varkappa$. The vector which is a basis
for one-dimensional space ${\rm ker}\,T_{\mu_1+1}$ gives coefficients of the first essential polynomial $Q_1(z)$ which is the minimal degree denominator of the Pade approximant.

Further we will denote by $T_{\mu_1+1}^{(k)}$ the matrix which is obtained from the matrix $T_{\mu_1+1}$ by deleting of the $k$th column, $k=1,\ldots,\varkappa+1$.

\begin{theorem} \label{vancoefQ}
Let $Q_1(z)=q_0+q_1z+\cdots+q_{\varkappa}z^{\varkappa}$ be the denominator with the minimal degree. Then
$q_k=0$ iff ${\rm rank}\,T_{\mu_1+1}^{(k+1)}=\varkappa-1,$ $k=0,\ldots, \varkappa.$
\end{theorem}
{\bf Proof.}
Let $q_k=0.$ We have ${\rm rank}\,T_{\mu_1+1}^{(k+1)}\leq \varkappa.$
Suppose ${\rm rank}\,T_{\mu_1+1}^{(k+1)}=\varkappa.$ Then $T_{\mu_1+1}^{(k+1)}$ has the trivial kernel, i.e. $q_0=q_1=\ldots=q_{k-1}=q_{k+1}=\ldots=q_{\varkappa}=0.$
Thus, $Q_1(z) \equiv 0$. Since it is impossible, we have ${\rm rank}\,T_{\mu_1+1}^{(k+1)}\leq \varkappa-1.$

Assume that ${\rm rank}\,T_{\mu_1+1}^{(k+1)}<\varkappa-1.$ Then $T_{\mu_1+1}^{(k+1)}$ has a multidimensional kernel, which can be embedded in the kernel of $T_{\mu_1+1}$ by a natural way.
But ${\rm ker}\, T_{\mu_1+1}$ is one-dimensional, hence our assumption is false and
${\rm rank}\,T_{\mu_1+1}^{(k+1)}=\varkappa-1.$

On the other hand, if ${\rm rank}\,T_{\mu_1+1}^{(k+1)}=\varkappa-1$, then
${\rm ker}\,T_{\mu_1+1}^{(k+1)}$ is one-dimensional and after the natural embedding in  one-dimensional space ${\rm ker}\,T_{\mu_1+1}$ we get $q_k=0.$

\halmos

Denote by $T_{\mu_1+1}^{[k]}$ the matrix which is obtained by inserting the row $(c_k\ c_{k-1}\ \ldots \ c_{k-\varkappa})$ at the beginning of the matrix $T_{\mu_1+1}$, $k=0,\ldots,m$.

\begin{theorem} \label{vancoefP}
Let $P_1(z)=p_0+p_1z+\cdots+p_{m}z^{m}$ be the numerator corresponding to the denominator $Q_1(z)$.
Then $p_k=0$ iff ${\rm rank}\,T_{\mu_1+1}^{[k]}=\varkappa,$ $k=0,\ldots, m.$
\end{theorem}
{\bf Proof.}
It follows from (\ref{systemnom}) that $p_k=c_k q_0+c_{k-1}q_1+\cdots+c_{k-\varkappa}q_{\varkappa},$ $k=0,\ldots,m.$ Hence,
$$
T_{\mu_1+1}^{[k]} \begin{pmatrix}
q_0 \\
q_1\\
\vdots\\
q_{\varkappa}
\end{pmatrix} =\begin{pmatrix}
p_k\\
0\\
\vdots\\
0
\end{pmatrix}.
$$

If $p_k=0$, then matrix $T_{\mu_1+1}^{[k]}$ has a nontrivial kernel, hence ${\rm rank}\,T_{\mu_1+1}^{[k]}<\varkappa+1.$ Since ${\rm rank}\,T_{\mu_1+1}=\varkappa$, we have ${\rm rank}\,T_{\mu_1+1}^{[k]}=\varkappa.$

On the other hand, if ${\rm rank}\,T_{\mu_1+1}^{[k]}=\varkappa$ then the inserted row
$(c_k\ c_{k-1}\ \ldots \ c_{k-\varkappa})$ is a linear combination of the rows of $T_{\mu_1+1}$. Therefore, $c_k q_0+c_{k-1}q_1+\cdots+c_{k-\varkappa}q_{\varkappa}=0,$  i.e. $p_k=0$.

\halmos

From Theorems \ref{vancoefQ}, \ref{vancoefP} we can obtain the multiplicities of $z=0$, $z=\infty$ as the roots of $Q_1(z),$ $P_1(z),$ i.e. the deficiency index
of $f_{m,n}(z)$ and the degrees of $Q_1(z),$ $P_1(z)$.

In particularly, the next result on the deficiency index is now evident.

\begin{corollary} \label{deficiencyindex}
The number $\delta$ is the deficiency index of $f_{m,n}(z)$ iff
$${\rm rank}\,T_{\mu_1+1}^{(1)}=\ldots={\rm rank}\,T_{\mu_1+1}^{(\delta-1)}=\varkappa-1,\ \  {\rm rank}\,T_{\mu_1+1}^{(\delta)}=\varkappa.$$
\end{corollary}

We would like to end this section with the following remark.
As is well known in the case of {\it normal} Pade approximant
$f_{m,n}(z)=\frac{P_{m,n}(z)}{Q_{m,n}(z)}$ the degrees of $P_{m,n}(z)$ and $Q_{m,n}(z)$ are equal to $m$ and $n$, respectively, and $\delta=0$.
Thus, in the normal case degrees of the numerator
and the denominator and the deficiency index $\delta$ are known.
Our modified Definition \ref{modpade} and the
results of this section allow us to determine
them not only in the normal case, but also in the non-normal one.

\section{Algorithm}\label{our}

In this section we present our algorithm for computing a reduced Pade
approximant with the minimal degree denominator.
As we have seen in the
previous section, in order to realize the algorithm,
we have to find the rank of matrix $T_m$, and
the null space of $T_{\mu_1+1}$. Moreover, in order to delete
vanishing coefficients of the denominator and the numerator we have to
find ranks of matrices $T_{\mu_1+1}^{(k)}$, $k=1,\ldots,\varkappa+1$,
and $T_{\mu_1+1}^{[k]}$, $k=0,\ldots,m$.
In practice, due to rounding and measuring errors and finite
computer precision, the elements of these matrices are perturbed with error,
so that we have to determine
the rank and the null space of the original matrix  from the
perturbed matrix.

For computer calculations we need to define the numerical
rank and the numerical null space of $A$.
Recall these concepts which are crucial for the algorithm.

The definition of the numerical rank was first given by Golub, Klema
and Stewart \cite{Golub}. We will use the simplified definition (as
in \cite{Foster}, \cite{Golub2}, p. 72). The
{\it numerical $\varepsilon$ rank} of an $M\times N$ matrix $A$ with
respect to the threshold $\varepsilon>0$ is defined as the smallest
rank of all matrices within a $2$-norm distance $\varepsilon$ of
$A$. Namely,
\[
{\rm rank}\,(A,\varepsilon)=\min_B \{{\rm rank}\, B: ||A-B||_2\leq
\varepsilon\}.
\]

The numerical $\varepsilon$ rank may be characterized in terms of
the singular value decomposition (SVD). According to Stewart
\cite{Stewart}, \flqq the singular value decomposition is the {\it
cr\`{e}me de la cr\`{e}me} among rank-revealing decomposition\frqq .

 Let us recall the definition of the SVD. Every $M\times N$ complex
matrix $A$ can be represented in the form
$$A=U\Sigma V^{H},$$ where $U$, $V$ are unitary matrices ($ ^H$
means the Hermite conjugation), $\Sigma$ is an $M\times N$ matrix in
which the upper $N\times N$ block is a diagonal matrix with all
entries real and sorted in descending order. The diagonal entries
$\sigma_1 \geq \ldots \geq \sigma_N$ of the matrix $\Sigma$ are
called the singular values.

 In the case of exact calculations, we have $\sigma_1 \geq
\ldots \geq \sigma_r>\sigma_{r+1}=\ldots=\sigma_N=0,$ where $r$ is
the rank of $A$. In practice, we have
\[
\sigma_1 \geq \ldots \geq \sigma_r
>\varepsilon \geq \sigma_{r+1} \geq \sigma_N\geq 0.
\]
Let us denote $A_r=U\Sigma_r V^H$ with $\Sigma_r={\rm
diag}\{\sigma_1,\ldots, \sigma_r,0,\ldots,0\}$; then
$||A-A_r||_2=\sigma_{r+1}$ and ${\rm rank}\,(A,\varepsilon)={\rm
rank}\, A_r=r$ \cite{Golub}. Moreover, $A_r$ is the nearest matrix
to $A$ (with respect to the 2-norm) with rank $r.$ Therefore the
null space of $A_r$ is called the numerical null space of $A$ within
$\varepsilon$. The null space of $A_r$ is spanned by
$\{v_{k+1},\ldots, v_{N}\}$, where $v_k$ is the $k$th column of the
matrix $V.$ The ratio $\gamma=\frac{\sigma_r}{\sigma_{r+1}}$ is
called the numerical rank gap. The  numerical rank of the matrix $A$
may be estimated reliably in the case of \flqq well defined
numerical rank\frqq \, that is when $A$ has a single well-determined
gap between large and small singular values.

Thus, in terms of singular value decomposition, numerical
$\varepsilon$ rank is the number of singular values greater than the
given threshold $\varepsilon.$ There is no uniform threshold for all
applications. The user must make a decision on the threshold
$\varepsilon$ based on the nature of the applications.

The {\tt rank(A,tol)} function in {\it Matlab} returns the number of singular values of $A$ that are larger than the threshold (tolerance) {\tt tol}.
The default tolerance is  $\max(M,N)*{\rm eps}({\rm norm}(A)).$
Here ${\rm norm}(A)$ indicates the Euclidean norm of $A$
and ${\rm eps}({\rm norm}(A))$ is approximately $2.2 \times
10^{-16}$ times ${\rm norm}(A).$ This choice is usually a good
choice if the errors in the matrix elements are due to computer
arithmetic and if there is a sufficiently large gap in the singular
values around this tolerance.

Although the SVD is the most widely used method for determination of
the numerical rank  and the numerical null space, there are
alternative methods like URV decomposition, LU decomposition or QR
factorization with column pivoting.

Now we can present our algorithm for computing $(m,n)$ reduced Pade approximant $f_{m,n}(z)$ for $f(z)$ at the point $z=a$.

{\bf Algorithm.}

{\it Initialization:}

\ \ \ \ \ \ \parbox[t]{0.9\textwidth}{
$f(z):=$ the approximated function;

$(m,n):=$ the order of the Pade approximant;

$a:=$ the center of expansion;

$d:=0$, if we want to delete vanishing coefficients of numerator and the denominator, else $d:=1$.
}

{\bf 1.}  Compute the Taylor coefficients $c_0,\ldots,c_{m+n}$ of
$f(z)$ at the point $z=a$.

{\bf 2.}  Form the Toeplitz matrix
$T_m=\|c_{i-j+m}\|_{\mbox{\scriptsize${\begin{array}{l}
i=1,\ldots,n+1,\\
j=1,\ldots,n.
\end{array}}$}}$

{\bf 3.} Determine its rank $\varkappa={\rm rank}\,T_{m}$ and the index
$\mu_1=m-n+\varkappa$.

{\bf 4.} Form matrix
$T_{\mu_1+1}=\|c_{i-j+\mu_1+1}\|_{\mbox{\scriptsize${\begin{array}{l}
i=1,\ldots,m+n-\mu_1,\\
j=1,\ldots,n-m+\mu_1+1.
\end{array}}$}}$

{\bf 5.} Find a basis for its one-dimensional kernel ${\rm ker}\,T_{\mu_1+1}$.

The obtained vector $(q_0, q_1, \ldots,  q_{n-m+\mu_1})^T$
is a vector formed from the coefficients of the minimal degree denominator
$Q_1(z)=\sum\limits_{k=0}^{n-m+\mu_1}q_k(z-a)^k$ of the Pade approximant.

{\bf 6.} If $d=0$, then
for $k=0,\ldots,\varkappa$, do:

\ \ \ \ \ \ \parbox[t]{0.9\textwidth}{
{\bf a.} Form the matrix $T_{\mu_1+1}^{(k+1)}$.

{\bf b.} Find its rank: ${\rm rank}\,T_{\mu_1+1}^{(k+1)}$.

{\bf c.} If ${\rm rank}\,T_{\mu_1+1}^{(k+1)}=\varkappa-1,$ then replace $q_k$ in $(q_0, q_1, \ldots,  q_{n-m+\mu_1})^T$
by zero.
}

\ \ \ Else {\bf step 6} should be omitted.

{\bf 7.} For $k=0,\ldots,m$ do $p_k=c_k q_0+c_{k-1}q_1+\cdots+c_{k-\varkappa}q_{\varkappa}$ and form the vector
$(p_0,\ldots,p_m)^T$ consisting of the coefficients of the numerator
$P_1(z)=\sum\limits_{k=0}^{m} p_k (z-a)^k$ corresponding to $Q_1(z).$

{\bf 8.}
If $d=0$, then
for $k=0,\ldots,m$, do:

\ \ \ \ \ \ \parbox[t]{0.9\textwidth}{
{\bf a.} Form the matrix $T_{\mu_1+1}^{[k]}$.

{\bf b.} Find its rank: ${\rm rank}\,T_{\mu_1+1}^{[k]}$.

{\bf c.} If ${\rm rank}\,T_{\mu_1+1}^{[k]}=\varkappa,$ then replace $p_k$ in $(p_0,\ldots,p_m)^T$
by zero.
}

\ \ \ \ Else {\bf step 8} should be omitted.

End of {\bf Algorithm.}

{\it Output:} $f_{m,n}(z)=\frac{P_1(z)}{Q_1(z)}$.

{\it Comments:}

1) {\bf Step 1} consisting of computing the
Taylor coefficients is very important for the algorithm.
It is well known (see, for example, \cite{Baker}) that accuracy
in the given coefficients $c_k$ is essential for Pade approximants.
Usually most computing effort goes into calculation of the coefficients
rather than Pade approximants, and so the coefficients
should be calculated with the greatest possible accuracy.
The way of their calculating depends
essentially on the problem under consideration.
For example, in the Pade -- Laplace method coefficients $c_k$
are found from experimental data by numerical integration. So the accuracy
in the coefficients depends not only on the chosen method of integration
but also on the quality of the experiment.

2) On {\bf step 3}  we find the essential
index $\mu_1$ according to Proposition \ref{indformula}. The first essential polynomial, which
is the minimal degree denominator of the Pade approximant, is found on {\bf step 5}
according to Theorem \ref{param}.
Theorems \ref{vancoefQ}, \ref{vancoefP} are used  on {\bf steps 6, 8}, respectively.
The output approximant $f_{m,n}(z)$
possesses properties listed in Theorem \ref{proper}.

We have implemented the algorithm in {\it Maple}
(procedure {\it ReducedPade(f,m,n,a,d)}) and now we would like to repeat
Examples \ref{Fr}, \ref{Fr1} from Section \ref{Prel}.

\begin{example}
Let us find by {\it ReducedPade} the diagonal $(4,4)$ Pade
approximant $f_{4,4}$ for $f(z)=\frac{(z+1)(z-2)}{(z+2.1)(z-1)}$ at the point
$z=0$ via the following commands.

f := $\frac{(z+1)\cdot(z-2)}{(z+2.1)\cdot(z-1)}:$

app := ReducedPade(f, 4, 4, 0, 0);

\begin{center}
$\frac{0.7773220744+0.3886610375z-0.3886610371z^2}{0.8161881781-0.427527140654618832z -0.388661037357416250z^2}$                                                             \end{center}

fsolve(numer(app), z, complex);
\vspace{-2mm}
\begin{center}
\footnotesize  $ -0.9999999998, 2.000000001$
\end{center}
\vspace{-2mm}

fsolve(denom(app), z, complex);
\vspace{-2mm}
\begin{center}
\footnotesize $-2.099999999, 1.000000000$
\end{center}

\end{example}

\begin{example}

Let us find the $(2,3)$ Pade approximant for the function
$f(z) = \frac{z+1.01}{(z+2)(z-2.01)}$ at the point $z=0.$

f := $\frac{z+1.01}{(z+2)\cdot(z-2.01)}:$

app := ReducedPade(f, 2, 3, 0, 0);

\begin{center}
$\frac{-0.2438127455 - 0.2413987580z}{0.9704230069 + 0.00241398781766979931z - 0.241398758054323787z^2}$                                                             \end{center}

fsolve(numer(app), z, complex);
\vspace{-2mm}
\begin{center}
\footnotesize  $ -1.010000000$
\end{center}
\vspace{-2mm}

fsolve(denom(app), z, complex);
\vspace{-2mm}
\begin{center}
\footnotesize $-1.999999999, 2.010000000$
\end{center}

\end{example}

As it can be seen,
the obtained approximants are identical to the approximated functions and Froissart doublets do not appear.
Note that in the given examples the value of the parameter $d$ was equal to zero.
Hence, the vanishing coefficients of the numerators and the  denominators were deleted. The next example shows that such deleting is desirable.

\begin{example}
Let us find the $(3,7)$ Pade approximant for the function
$f(z) = \frac{z+1.01}{z^4+3z^2-4.01}$ at the point $z=0.$

f := $\frac{z+1.01}{z^4+3z^2-4.01}:$

app := ReducedPade(f, 3, 7, 0, 0);

\begin{center}
$\frac{0.1977728838 + 0.1958147364z}{-0.7852170931 + 0.587444214531195330z^2+ 0.195814737672938278z^4}$
\end{center}

fsolve(numer(app), z, complex);
\vspace{-2mm}
\begin{center}
\footnotesize  $ -1.010000000$
\end{center}
\vspace{-2mm}

fsolve(denom(app), z, complex);
\vspace{-2mm}
\begin{center}
\footnotesize  $ -1.000999098, -2.000499738 I, 2.000499738 I, 1.000999098$
\end{center}
\vspace{-2mm}

app := ReducedPade(f, 3, 7, 0, 1);

\begin{center}
\footnotesize $(0.1977728838+0.1958147377z-1.48862691\, 10^{-11}z^3)\slash
(-0.7852170931- 5.23710724786852211 \, 10^{-9}z + 0.587444214531194664 z^2 - 1.13198619922094679 \, 10^{-9}  z^3  + 0.195814737672938056 z^4- 1.67404646568558580\, 10^{-10}z^5)$
\end{center}

fsolve(numer(app), z, complex);
\vspace{-2mm}
\begin{center}
\footnotesize  $ -1.146906047\, 10^{5} , -1.009999994, 1.146916147\, 10^5$
\end{center}
\vspace{-2mm}

fsolve(denom(app), z, complex);
\vspace{-2mm}
\begin{center}
\footnotesize  $  -1.000999095, -1.728850777\, 10^{-9} - 2.000499738 I, -1.728850777 \,10^{-9} + 2.000499738 I, $ $ 1.000999101, 1.169709095\, 10^9$
\end{center}
\vspace{-2mm}
\end{example}

\section{Conclusions}

This work was motivated by our intention of avoiding the Froissart doublets appearance
in signal processing by Pade -- Laplace method. To do this, we have modified the classical
definition of Pade approximant by adding the requirement of the minimal degree of its denominator.
It turned out that this new definition of {\it reduced Pade approximant} allows us to avoid Froissart
doublets induced by computer roundoff.

The reduced Pade approximant can be easily obtained by the proposed algorithm
and always has nice properties which classical Pade
approximant possesses only in the normal case: the denominator is unique up to a constant factor,
the numerator and the denominator have not common non-zero roots,
their degrees and the deficiency index are known exactly.

\end{document}